\documentclass{article}
\usepackage[utf8]{inputenc}
\usepackage{amsmath}
\usepackage{amsfonts}
\usepackage[margin=0.5in]{geometry}

\usepackage{graphicx}
\usepackage{tabu}

\title{Bounds for Combinatorial Types of Non-Attacking Riders}
\author{Grant Jensen}
\date{22 June 2020}

\begin{document}

\maketitle
\section{Introduction/Background}
\large{\null\quad Consider first an n by n infinite chessboard upon which is placed a queen. Now, if a second queen is placed, we dictate that these queens cannot be in line of fire of each other. With this assumption, it is clear that the second piece can be placed in 8 different locations with respect to the first piece; however, it will be explained later that this only results in 4 combinatorial types. We now expand this concept to an infinite plane upon which "pieces" or "riders" may be placed with no grid present. These two terms are interchangeable. These riders will have some r number of moves denoting the number of different straight line movements that can be made. From the example above, the queen would have 4 moves, and a rook would have 2. Pawns, knights, and kings would not exist in this scenario because they do not have straight line movements, and cannot continue in a direction indefinitely. One can also imagine a piece with 5+ moves. Such a piece could simply move in 5 or more straight directions indefinitely. Also note that these moves need not be set at a uniform angle apart. We wish to find a formula which will output the number of combinatorial types for any given number of riders and moves. Note that we will not mix riders with different movement patterns. For example we will not solve the number of combinatorial types for a rider with 4 moves placed on a plane with a rider with 2 moves.}
\begin{figure}[h]
\centering
\includegraphics[width=5cm, height=5cm]{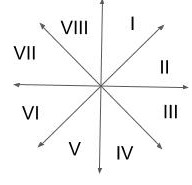}
\caption{Labeling the 8 regions for a piece with 4 moves.}
\end{figure}

\section{General Syntax}
\large{\quad Two non-attacking pieces are said to have the same combinatorial type if for each pair $\mathbb{P}_i$ and $ \mathbb{P}_j$, $\mathbb{P}_j$ lies in the same region of the board with respect to $\mathbb{P}_i$ in both configurations (Hanusa). We will first introduce some notation. Without loss of generality, let every piece contain a move which is vertical. Then, denote region I as the region directly clockwise to the the upper half of this vertical move. Let region II be the region clockwise of region I. Continue clockwise in this manner, labeling every region of the piece. Given a combinatorial type, (such as Figure 1) we record it by first choosing a piece, call it $\mathbb{P}_1$. Then, choose any other piece, $\mathbb{P}_2$, and record the region of  $\mathbb{P}_1$ that $ \mathbb{P}_2$ is in. Begin a list, say x. Let the first index of x be this number we recorded, ie x=(1). Then, choose a third piece, $\mathbb{P}_3$, and let the second index of x be the region of $\mathbb{P}_1, \mathbb{P}_3$ is in. Let the third index of x be the region of $\mathbb{P}_2, \mathbb{P}_3$ is in, ie x=(1,6,5). Continue this process, so for some $\mathbb{P}_i$, it is recorded in the indices $[\frac{i(i-1)}{2}+1,\frac{i(i+1)}{2}]$ of x, the j-th of which being its location relative to $\mathbb{P}_j$, $j<i$.  Note that one combinatorial type can be recorded in multiple different ways.\\ 

\begin{figure}[h]
\centering
\includegraphics[width=16cm, height=8cm]{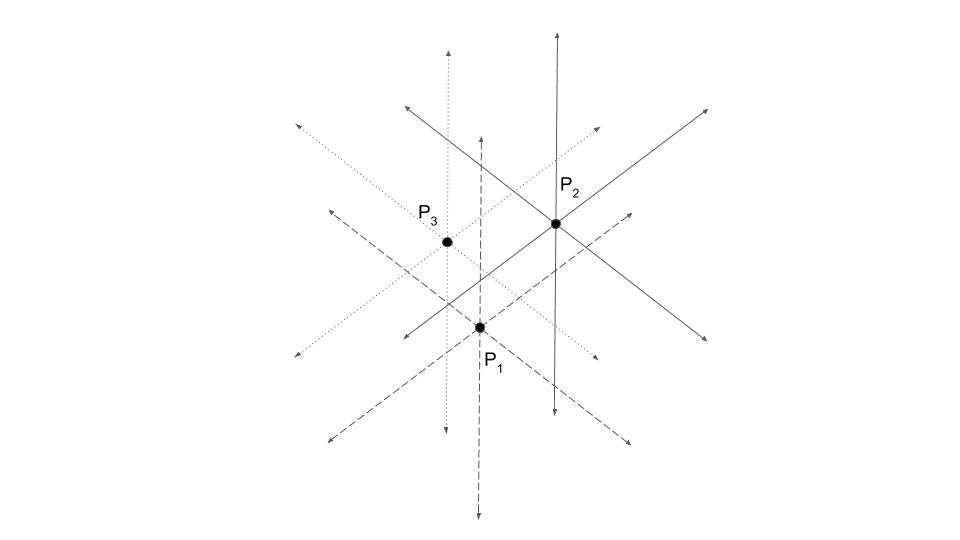}
\caption{This would typically be recorded as (1,6,5). However, it could also be (6,1,2) if we label in the order ($\mathbb{P}_1,\mathbb{P}_3,\mathbb{P}_2$).There are 8 different ways to order 3 pieces, thus there are 8 possible ways to record the above figure. We say there are 8 permutations of the above combination. }
\end{figure}

\textbf{Lemma 1}: There are q! ways to record a combinatorial type with q riders. Call each recording of the combinatorial type a permutation.\\ \\
\null\quad Proof: This is trivial. In order to record a different permutation of the same combination, simply record the pieces in a different order. In other words, choose a different piece to denote as $\mathbb{P}_i$ when recording the permutation. There are q! ways to order the list of integers from one to q, so there must be q! different permutations of one combinatorial type. All that remains is to show that each permutation must be unique. Suppose not two orderings of the pieces are unique, but their permutation is the same. Ie one labels pieces in the order (1,2,3) and results in the permutation (1,6,5), while the other has pieces labeled as (2,1,3) but also results in the permutation (1,6,5). This is clearly impossible, thus we have proved our claim.\\ \\
\null\quad Therefore to find the number of combinatorial types for q riders and r moves, it suffices to find the number of permutations.}

\section{Lower Bound}
\large{Define spaces as the largest areas for which no moves intersect. Consider Figure 2 which has 34 spaces. \\ \\
\null\quad \textbf{Lemma 2}: Every combinatorial type of q riders and r moves has the same number of spaces.\\ 
Proof: A space is created by the intersection of a space with a move. Proof by induction on the number of riders, q. Base case: If there is one rider, then since there is only one combinatorial type, the number of spaces is the same for every type. Inductive step: Suppose for r moves, every combinatorial type with q-1 riders has the same number of spaces. Regardless of where this new q-th rider is placed, each of its moves will intersect with (r-1)(q-1) moves from the other q-1 riders. (Each move will not intersect the parallel move of the other riders). Note that the q-th rider will transform the space it is placed in into (2r-1) spaces. Also, its moves will create r(r-1)(q-1) additional spaces since r moves intersect (r-1)(q-1) lines. So, we can create a formula for the number of spaces. Define s(q,r) represent the number of spaces available for q riders and r moves. Then, s(q,r)=(2r-1)+r(r-1)(q-1)+s(q-1,r). By the inductive step, each combinatorial type has the same number of spaces.\\ \\
\null\quad Moreover, we now have an inductive formula for the number of spaces. This can be re-written as $s(q,r)=\frac{q(q+1)}{2}(r^2-r)+q(-r^2+3r-1)+1$.Next, define p(q,r) as the number of permutations for q riders and r moves. It is clear that $p(q,r)= \prod\limits_{n=1}^{q-1}s(n,r)$ via Lemma 2. Thus, conclude via Lemma 1 that \[t(q,r)=\frac{1}{q!}\cdot \prod\limits_{n=1}^{q-1}[\frac{n(n+1)}{2}(r^2-r)+n(-r^2+3r-1)+1]\]
}
\section{Upper Bound}
\large{\null\quad While this agrees with previous formulas (Hanusa), (Kot$\check{e}\check{s}$ovec) for low q and r, it acts as a lower bound once we reach $q\geq 4, r\geq 3$. This is due to a slight error within Lemma 2. Suppose q-2 riders with r moves have already been placed on the board. When the q-1 rider is placed on the board in certain spaces, two different sets of spaces can be created based upon where in the space the q-1 rider is placed. Thus, when the final q-th rider is placed, it may have a larger than normal set of spaces to be placed in depending on where the q-1 rider is placed. Figure 3 demonstrates this.  \\

\begin{figure}[h]
\centering
\includegraphics[width=16cm, height=8cm]{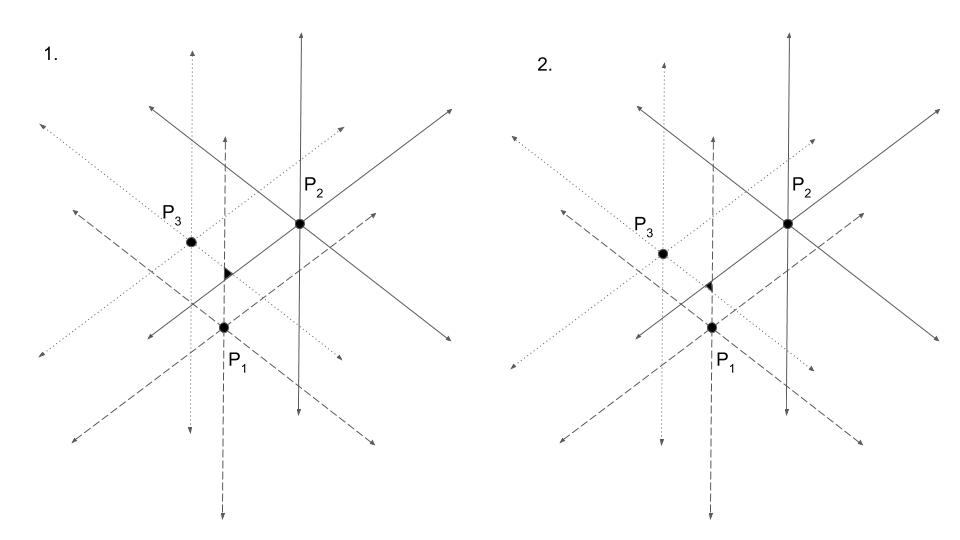}
\caption{First note the in both diagrams we have the same recording, that is (1,6,5). In diagram 1, we have a blackened triangle. If a 4-th rider were placed in this triangle, we would record it as (1,6,5,1,5,3). This location does not exist in diagram 2. Moreover, although the first three riders are combinatorically equivalent in both diagrams, we record diagram 2 as (1,6,5,6,4,2). Once again, this space does not exist in diagram 1. Thus s(4,3) will be larger than expected due to extra spaces with the same orientation of the first three pieces.}
\end{figure}

Determining how many of these "problem spaces" occur given some q and r is a difficult problem, and is yet to be solved. However, an upper bound can be created. These problem spaces only occur at intersections of three or more moves with less than q riders placed, hence the bounds of $q\geq 4, r\geq 3$. (This part is clear, since Lemma 2 accounts for all intersections of 2 lines). Thus, if we simply add to s(q,r) every instance of three move intersections, we can find a formula for t(q,r) given any q and r. This however, is also very difficult, thus we will create an upper bound for the number of 3 line intersections with some simple geometry. \\ \\
\null\quad Note that for q-2 riders and r moves, there are $r(r-1)\frac{(q-2)(q-3)}{2}$ intersections of 2 moves. This can be seen by looking at one rider at a time. First, count the number of intersections created by the q-2 rider. Each of its r moves will intersect with r-1 moves from the other q-3 riders. Next, the q-3 rider will have r moves intersect with r-1 moves from the other q-4 riders (not the q-2 rider). Summing these up yields $=r(r-1)[(q-3)+(q-4)+...+1]=r(r-1)\frac{(q-2)(q-3)}{2}$. Next, assume that the q-1 rider can reach any 2 move intersection with any move which is not parallel to either of the two moves forming the intersection. (This is why this will be an upper bound). So, when placing the q-1 rider, each move can interact with $\frac{r-2}{r}\cdot r(r-1)\frac{(q-2)(q-3)}{2}$ 2-intersections. If we define a(q,r) as the number of problem spaces for q moves and r riders, we conclude that $a(q,r)=r(r-1)(r-2)\frac{(q-2)(q-3)}{2}$. This results in the following formula:
\begin{gather}
    \frac{1}{q!}\cdot \prod\limits_{n=1}^{q-1}[\frac{n(n+1)}{2}(r^2-r)+n(-r^2+3r-1)+1]\leq t(q,r)\leq  \\ \frac{1}{q!}\cdot \prod\limits_{n=1}^{q-1}[\frac{n(n+1)}{2}(r^2-r)+n(-r^2+3r-1)+1+r(r-1)(r-2)\frac{(n-2)(n-3)}{2}] 
\end{gather}

\large{Where (2) only applies when $q\geq 4$ and $r\geq 3$. One might ask what occurs at intersections of 4 or more lines. As seen in Figure 3, an intersection of 3 lines creates one problem space. It can be seen clearly that an intersection of x lines (where x$\geq$3) creates (x-2) problem spaces. However, a(q,r) counts an intersection of x lines as ${x-1 \choose 2}$ problem spaces. (There are ${x-1 \choose 2}$ 2-intersections in an (x-1)-intersection, then the q-th rider makes it a 3-intersection which equals 1 problem space). Clearly ${x-1 \choose 2}\geq (x-2)$ for $x\geq 2$, so we need not worry about large intersections since they will still be counted below the upper bound. This concludes the proof. Below are the results of these equations compared to Hanusa's results.}
\begin{table}[h]
\centering
\begin{tabu}{|c|[2pt]c|c|c|c|c|c|}
    \hline
     q$\backslash$ r&1&2&3&4&5&6  \\ \tabucline[2pt]{-}
     1&1&1&1&1&1&1 \\ \hline
     2&1&2&3&4&5&6 \\ \hline
     3&1&6&17&36&65&106 \\ \hline
     4&1&24&144.5$\leq$ 151$\leq$ 289&522$\leq$ 574$_\mathbb{Q}\leq$ 2088 &1430$\leq$?$\leq$ 10010&3286$\leq$?$\leq$36146\\ \hline
     5&1&120&1647.3$\leq$1899&10544.4$\leq$14206$_\mathbb{Q}$&44902$\leq$?&147870$\leq$?  \\ 
     ~&~&~&$\leq$3641.4&$\leq$52200&$\leq$434434&$\leq$2494074\\ \hline
     6&1&720&23611.3$\leq$31709&274154.4$\leq$501552$_\mathbb{Q}$& 1840982$\leq$?&8773620$\leq$?\\ 
     ~&~&~&$\leq$ 63117.6&$\leq$1983600&$\leq$ 30844814&$\leq$297626164 \\ \hline
\end{tabu}
\caption{For q$\leq$3 and r$\leq $2, the upper bound formula agrees exactly with results previously found by Kot$\check{e}\check{s}$ovec, and equations from Hanusa. For the other entries, the middle value is from Kot$\check{e}\check{s}$ovec's empirical formulas while the bounds come from the formulas above. $\mathbb{Q}$ denotes that the empirical formula only applies for riders with moves identical to a chess queen. Riders with 4 moves in different orientations may result in different values. }
\end{table}
\\ \\
\section{Conclusion}
\large{\null\quad Clearly much progress can be made on improving these bounds(especially the upper bound). The question that remains is how to calculate the number of these problem spaces given q and r. In addition, it is possible that different types of pieces (pieces with the same number of moves, but whose moves are a different angle apart) could result in a different number of combinatorial types despite having the same number of moves and riders!}
}

\section{References}
\large{
Hanusa, Christopher R. H., and Thomas Zaslavsky. “A q-Queens Problem. VII. Combinatorial Types of Nonattacking Chess Riders.” ArXiv.org, 11 June 2020, arxiv.org/abs/1906.08981.
\\ \\
V. Kot$\check{e}\check{s}$ovec, Non-attacking chess pieces (chess and mathematics) [\textit{$\check{S}$ach a matematika - po$\check{c}$ty rozm\'ist$\check{e}$n\'i neohro$\check{z}$uj\'ic\'ich se kamen$\overset{\circ}{u}$}]. Self-published online book, Apr. 2010; 6th ed., Feb. 2013. \\
http://www.kotesovec.cz/math.htm
}

\end{document}